\documentclass[11pt]{amsart}

\usepackage{lineno}
\usepackage{palatino}
\usepackage[all]{xy,xypic}
\usepackage{amsfonts,amsmath,oldgerm,amssymb,amscd,tikz-cd,enumerate}
\usepackage[a4paper,left=12em,right=12em,top=10em,bottom=10em]{geometry}

\makeatletter
\@namedef{subjclassname@2020}{%
	\textup{2020} Mathematics Subject Classification}
\makeatother

\newcommand{\ra}{\rightarrow}		

\newcommand{\by}[1]{\stackrel{#1}{\ra}}

\newcommand{\remove}[1]{}

\newcommand{\surj}{\ra\!\!\!\ra}	

\newcommand{\ol}{\overline}		

\newcommand{\iso}{\by \sim}

\newtheorem{theorem}{Theorem}[section]
\newtheorem{proposition}[theorem]{Proposition}
\newtheorem{lemma}[theorem]{Lemma}
\newtheorem{definition}[theorem]{Definition}
\newtheorem{corollary}[theorem]{Corollary}
\newtheorem{conjecture}[theorem]{Conjecture}
\newtheorem{example}[theorem]{Example}

\newcommand{\CO}{\mbox{$\mathcal O$}}	\newcommand{\CP}{\mbox{$\mathcal P$}}
	
\newcommand{\CS}{\mbox{$\mathcal S$}}

\newcommand{\MI}{\mbox{$\mathfrak I$}}

	\newcommand{\MT}{\mbox{$\mathfrak T$}}

\newcommand{\ma}{\mbox{$\mathfrak a$}}

	\newcommand{\p}{\mbox{$\mathfrak p$}}
\newcommand{\mq}{\mbox{$\mathfrak q$}}

\newcommand{\Spec}{\text{Spec}}	
\newcommand{\hh}{\text{ht}}
\newcommand{\rank}{\text{rank}}

\newcommand{\Aut}{\mbox{\rm Aut\,}} 
\newcommand{\Hom}{\text{Hom}\,}

\newcommand{\Um}{\textnormal{Um}}		

\newcommand{\E}{\text{E}}	
\newcommand{\GL}{\text{GL}}

\newcommand{\bp}{\begin{proposition}}
	\newcommand{\ep}{\end{proposition}}
\newcommand{\bl}{\begin{lemma}}
	\newcommand{\el}{\end{lemma}}
\newcommand{\bt}{\begin{theorem}}
	\newcommand{\et}{\end{theorem}}
\newcommand{\bc}{\begin{corollary}}
	\newcommand{\ec}{\end{corollary}}
\newcommand{\bd}{\begin{definition}}
	\newcommand{\ed}{\end{definition}}
\newcommand{\bco}{\begin{conjecture}}
	\newcommand{\eco}{\end{conjecture}}
\newcommand{\bma}{\begin{bmatrix}}
	\newcommand{\ema}{\end{bmatrix}}

\def\proof{\paragraph{Proof}}
\def\example{\refstepcounter{theorem}\paragraph{{\bf Example} \thetheorem}}
\def\notation{\paragraph{\bf Notation}}

\title [Subrings of polynomial rings and the conjectures of Eisenbud and Evans]{Subrings of polynomial rings and \\ the conjectures of Eisenbud and Evans}
\author{Sourjya Banerjee}
\address{Department of Mathematics and Statistics, Indian Institute of Science Education and Research Kolkata, Campus Road, Mohanpur, West Bengal 741246, India}
\email{sourjya.pdf@iiserkol.ac.in, sourjya91@gmail.com}
\date{\today}


\usepackage{hyperref}

\begin{document}
	\maketitle

\begin{abstract}
Let $R$ be a commutative Noetherian ring of dimension $d$. In 1973, Eisenbud and Evans proposed three conjectures on the polynomial ring $R[T]$. These conjectures were settled in the affirmative by Sathaye, Mohan Kumar and Plumstead. One of the primary objectives of this article is to investigate the validity of these conjectures over Noetherian subrings of $R[T]$ of dimension $d+1$, containing $R$. We formulate a class of such rings, which includes polynomial rings, Rees algebras, Rees-like algebras and Noetherian symbolic Rees algebras, and exhibit that all three conjectures hold for rings belonging to this class.

\end{abstract}

\section{Introduction}
Let $R$ be a commutative Noetherian ring of finite (Krull) dimension $d$. As the title of the article suggests, we commence by recalling some old conjectures on $R[T]$ proposed by Eisenbud and Evans \cite{EEC}, which are now theorems. For the definitions related to the following conjectures, one may refer to Definition \ref{def}.
\begin{enumerate}[ 1.]
	\item\label{c1} Let $M$ be a finitely generated $R[T]$-module such that $\mu_{\p}(M)\ge \dim(R[T])$ for all $\p\in \Spec(R[T])$, where $\mu_{\p}(M)$ is the minimal number of generators of $M_{\p}$. Then $M$ has a basic element.
	\item\label{c2} Let $P$ be a finitely generated projective $R[T]$-module of rank $\ge \dim(R[T])$. If $Q$ is another $R[T]$-module such that $R[T]\oplus P\cong R[T]\oplus Q$, then $P\cong Q$. In other words, the module $P$ is cancellative.
	\item\label{c3} Let $M$ be a finitely generated $R[T]$-module. Then $M$ can be generated by $e(M)$ many elements, where
\end{enumerate}

 {\setlength{\abovedisplayskip}{-20pt}
	\setlength{\belowdisplayskip}{4pt}$$\quad e(M)=\sup\{\mu_{\p}(M)+\dim(R[T]/\p):\p\in \Spec(R[T]) \text{ such that } \dim(R[T]/\p)\le d\}.$$}
Sathaye \cite{AMS} proved conjecture \ref{c3} for affine domains over infinite fields. Mohan Kumar \cite{NMK} settled conjecture \ref{c3}, in general. Plumstead settled the other two conjectures in \cite{P}, which is one of the primary interests of this article.

 While the aforementioned conjectures hold for the ring $R[T]$, replacing it with $R$ renders them invalid in general. For instance, if $R$ is the coordinate ring of the sphere of dimension $2$ over $\mathbb{R}$, then all three conjectures have negative answers. However, the validity of these conjectures on Noetherian rings of dimension $d+1$ lying between $R$ and $R[T]$ remains unknown to date. This motivates us to investigate these questions on such rings, which constitutes one of the main foci of this article.

We formulate a class of Noetherian subrings of $R[T]$ of dimension $d+1$, containing $R$, which are referred to as ``geometric subrings" [see Definition \ref{gsr}] and prove that all three conjectures are valid on rings belonging to this class. Our model examples of such rings are coming from Rees algebras and Ress-like algebras. The latter one typically appears in \cite{E-G} in order to provide counterexamples of the Regularity Conjecture by Eisenbud and Goto. Another kind of rings which are included in our class is the Noetherian symbolic Rees algebras.

\subsection{Literature survey}
Here we pause to reflect on the significance of the articles \cite{AMS}, \cite{NMK} and \cite{P} in the literature. Indeed, the impact of these theorems goes beyond the results themselves. The methods and techniques used to prove these results have paved the way for further development in the literature. In \cite{AMS}, one of the noteworthy steps in proving conjecture \ref{c3} for affine domains over an infinite field was the reduction of the conjecture to a question on ideals. This reduction has been further refined and utilized in \cite{NMK} to establish conjecture \ref{c3} in its entirety. Additionally, \textit{Sathaye's theorem} has played a crucial role in \cite[Section 4]{M} by enabling Murthy to attach a finitely generated module over smooth affine varieties on an algebraically closed field with a certain ``Segre class" in the Chow group of zero cycles. Furthermore, it has been proven that the module is generated by the number of elements estimated by Eisenbud-Evans if and only if the assigned Segre class vanishes.

	Recall that, an ideal $K\subset R$ is said to be \textit{efficiently generated} if $\mu(K/K^2)=\mu(K)$, where $\mu(-)$ stands for the minimal number of generators. One of the crucial steps of Mohan Kumar's proof \cite[Theorem 2]{NMK} of conjecture \ref{c3} is the following: let $I\subset R[T]$ be an ideal such that (1) $\hh(I)\ge 1$ and (2) $\mu(I/I^2)=d+1$, then $I$ is efficiently generated. To prove this, he used some deep-sheaf patching techniques along the line of Quillen's patching theorem \cite{Q}. In the same article, Mohan Kumar also provided a partial solution \cite[Theorem 5]{NMK} to another conjecture posed by Murthy \cite{MCI} on complete intersections. However, Mohan Kumar's proof of the aforementioned step contains one of the key ideas, that he crucially used to provide a partial solution to Murthy's conjecture. Mohan Kumar's theorems \cite[Theorem 2 and 5]{NMK} offer a fundamental technique for tackling such problems, which has been further developed by Mandal \cite{SM1}. Mandal employed a refined version of Mohan Kumar's methodology in \cite{SM} to solve another question due to Nori, which was later reworked by Mandal and Sridharan in \cite{MR}. It was the insightful work of Bhatwadekar and Sridharan in \cite{SMBB3} that established a ``subtraction principle" in which \cite{MR} played a pivotal role. This subtraction principle serves as one of the key building blocks for the development of the ``Euler class theory", which acts as an obstruction group for the \textit{splitting problem} of projective modules of \textit{top} rank. The only improvement of \textit{Mohan Kumar's bound} to date has been achieved recently by Das \cite{MKD}.

On the other hand, Plumstead's proof of conjecture \ref{c2} has had a significant impact on the literature on its own. In \cite{P}, he solved the first two conjectures in the affirmative and also gave an independent proof of the third. Notably, he used the cancellation result \cite[Theorem 1]{P}, to prove the other two theorems. One of the most remarkable ideas of \cite{P} was the introduction of the concept of generalized dimension functions, which has had a profound impact on tackling many problems in the literature. Another key aspect of Plumstead's work is the elegant implementation of deep sheaf patching techniques, which are closely related to one of Quillen's key ideas in his groundbreaking article \cite{Q} on Serre's conjecture. Later, Mandal \cite{m82} extended Plumstead's result over Laurent polynomial rings in one variable. This was further extended by Bhatwadekar-Roy \cite{SMBAR} and Rao \cite{RROR} on some overrings of polynomial rings.
\subsection{Our approach} In our approach to handling geometric subrings $A\subset R[T]$, we follow Plumstead's footsteps of obtaining two covers of $\Spec(A)$. On the \textit{critical} cover, we establish the results with the help of (1) Eisenbus-Evans's foundational theorem on general stability arguments \cite{EE} and (2) Plumstead's concept of generalized dimension function \cite{P}. Finally, we apply sheaf patching techniques to establish our desired results on $\Spec(A)$. Notably, we observe that Plumstead's results on patching \cite[Lemma 2, Proposition 1 and 2]{P} hold beyond the polynomial rings. To achieve this, we get back to Quillen's initial approach as described in \cite{Q}. When dealing with the cancellation problem in Proposition \ref{canl}, we employ Plumstead's patching lemma \cite[Lemma 1]{P} to conclude our proof. However, while Plumstead's approach leading up to the patching in this proposition serves as a foundation, it requires some additional support. One of the true motivations of this came from \cite{m82}.

\subsection{Main results} The article is organized as follows: Section \ref{2} covers basic definitions and preliminary results necessary for proving the main theorems. In Section \ref{3}, the main theorem is Theorem \ref{ebe}, where we establish the existence of basic elements in certain modules, as indicated in conjecture \ref{c1}, over a geometric subring of $R[T]$.  Section \ref{4} focuses on the cancellation properties of projective modules. The main theorem of this section is Theorem \ref{cane}, which states that any projective module over a geometric subring of $R[T,f^{n}]$, of rank $\ge d+1$, is cancellative, where $f\in R[T]$ and $n\in \mathbb{Z}$. In Section \ref{5}, we establish that any finitely generated module over a geometric subring of $R[T]$ can be generated \textit{efficiently} in a view towards the estimation suggested by Eisenbud-Evans.

%Section \ref{6} is devoted in studying graded subrings $B$ of $R[T]$ containing $R$ of dimension $d+1$. Here we develop an analogy of Quillen-Suslin's theory of ``Local Global Principle" for the elementary subgroup $\E_n(B)$ of $\GL_n(B)$. As an application, we strengthen the results obtained in Section \ref{4}, specifically, we prove in Theorem \ref{ec} that any unimodular row in $B[X_1,...,X_m]$ of length $d+2$, can be completed to the first row of an elementary matrix. In Theorem \ref{cangr} we prove that any projective $B$-module (with trivial determinant) of rank $d+1$ is cancellative. In Theorem \ref{egi} we prove that if $I\subset A$ is an ideal such that $\mu(I/I^2)=\hh(I)=d+1$, then $\mu(I)=d+1$. As a consequence we show that every projective $B$-module of rank $d+1$ splits off a free summand of rank one. Additionally, in Theorem \ref{ist}, we prove that the injective stability of $\text{K}_1 (A[X_1,...,X_m])$ is $d+2$, which improves the existing results of Vaser{\v{s}}te{\u{\i}}n \cite{Va1} and Suslin \cite{AASSC}. 
%In section \ref{7} we demonstrate some applications of the main theorems, related to (1) an old question due to Roitman on a \textit{monic inversion principle}, and (2)  \textit{projective equivalence of ideals}.

\section{Preliminaries}\label{2}

This section summarizes several results and definitions from the literature that are used frequently in the article to prove the main theorems. We may restate or improve these results as required. Before further proceeding, in order to prevent any ambiguity, we introduce a set of conventions that will remain fixed throughout the entire article.
\smallskip

\paragraph{\bf{Convention}} The symbols $\mathbb{Z}$ and $\mathbb{N}$ will denote the set of all integers and non-negative integers respectively. In particular, we are assuming $0\in \mathbb{N}$. All rings considered in this article are assumed to be commutative Noetherian with $1\neq 0$. Unless otherwise stated, the symbol $R$ will always denote a commutative Noetherian ring of finite (Krull) dimension $\mathbf{d\ge 1}$. Any module considered in this article is assumed to be finitely generated.

We now turn our attention to the class of subrings of $R[T]$, which are one of the primary interests in this article. This class was introduced in \cite{BBP} and referred to as ``geometric subring". Hence we adopt their terminology to avoid further renaming.
\smallskip 

\begin{definition}\label{gsr}
Let $f\in R[T]\setminus R$ be a non-zero divisor and $n\in \mathbb{Z}$. A Noetherian subring $A$ of $R[T,f^n]$ containing $R$, is said to be a \textit{geometric subring} of $R[T,f^n]$, if (1) there exists a non-zero divisor $s\in R$ such that $A_s=R_s[T,f^n]$ and (2) $\dim(A)=d+1$.
\end{definition}

%		\bd
%		Let $A\subset B$ be two rings. The conductor ideal of $A$ with respect to $B$ is the ideal $\mc(B/A):=\text{Ann}_A(B/A)=\{a\in A:aA\subset B\}$.
%		
%		\ed

\smallskip

\example\label{ex1} Let us assume that $R$ is reduced and consider an ideal $I \subset R$ with $\hh(I)\ge 1$. Some interesting examples are the following.
\begin{enumerate}
	\item Rees algebra $R[IT]=\bigoplus_{n\in \mathbb{N} } I^nT^n$;
	\item Rees-like algebra $R[IT,T^2]= R\oplus IT\oplus RT^2\oplus I^3T^3\oplus RT^4\oplus \cdots $;
	\item Extended Rees algebra $R[IT,T^{-1}]=\bigoplus_{n\in \mathbb{Z} } I^nT^n$;
	\item Noetherian symbolic Rees algebra $R_S[IT]=\bigoplus_{n\in \mathbb{N}} I^{(n)}T^n$, where $I^{(n)}$ is the $n$-th symbolic power of $I$.
	\item Consider the ring $A=\mathbb{Z}[T^2+T^3,2T]$. Then we note that $A$ is a geometric subring of $\mathbb{Z}[T]$. Moreover, it is not difficult to check that $T^2\not\in A$. Hence $A$ is a geometric subring of $\mathbb{Z}[T]$ but not a graded subring of $\mathbb{Z}[T]$.
%	\item Let $h(T)\in R[T]$. We define $R[IT,h]:=\{f(T)+g(h(T)):f(T)\in R[IT],\text{ and } g(T)\in R[T]\}$.
\end{enumerate}

Before proceeding further, we revisit several definitions from the literature.
\smallskip

\bd\label{def} Let $M$ be a $R$-module.
\begin{enumerate}  
	\item  An element $x\in M$ is said to be a \textit{basic element} of $M$ at a prime ideal $\p\in \Spec(R)$ if $x\not \in \p M_{\p}.$ Let $\CS\subset \Spec(R)$. We call $x$ is a basic element of $M$ on $\CS$ if it is a basic element of $M$ at each prime ideal $\p\in \CS$. Whenever $\CS=\Spec(A)$, we omit mentioning $\CS$ explicitly.
	\item We denote the Eisenbud-Evans's estimation on the number of generators for $M$ as a module over $R$, by {\setlength{\abovedisplayskip}{1pt}
		\setlength{\belowdisplayskip}{3pt}\begin{equation*}	\,\,\,\,\,\,\,\,\,\,\,\,\,\,\,\,\,\,\,\,\,\, e(M):=\sup\{\mu_{\p}(M)+\dim(R/\p):\p\in \Spec(R)\text{ such that } \dim(R/\p)< \dim(R)\}.
	\end{equation*}}
	\item Let $\CS\subset \Spec(R)$ and let $d:\CS\to \mathbb{N}$ be a function. For two prime ideals $\p,\mq\in \CS$, we define a partial order $\p<<\mq$ if and only if $\p\subset \mq$ and $d(\p)>d(\mq)$.  We say that $d$ is a \textit{generalized dimension function} if for any ideal $I\subset R$, the set $V(I)\cap \CS$ has only finitely many minimal elements with respect to the partial ordering $<<$.
		\item The \textit{order ideal} of an element $m\in M$ is defined by $\CO_M(m):=\{\alpha(m)\in R:\alpha\in M^*=\Hom_R(M,R)\}.$ An element $m\in M$ is said a \textit{unimodular element} if $\CO_M(m)=R$. The set of all unimodular elements of $M$ is denoted by $\Um(M)$. If $M=R^n$, then we write $\Um_n(R)$ instead of $\Um(R^n)$.
	\item  Let $P$ be a projective $R$-module such that either $P$ or $P^*$ has a unimodular element. We choose $\phi\in P^*$ and $p \in P$ such that $\phi(p)=0$. We define an endomorphism $\phi_p$ as the composite $\phi_p:P\to R\to P$, where $R\to P$ is the map sending $1\to p.$ Then by a \textit{transvection} we mean an automorphism of $P$, of the form $1+\phi_p$, where either $\phi\in \Um(P^*)$ or $p\in \Um(P)$. By $\E(P)$  we denote  the subgroup of $\Aut(P)$ generated by all transvections.
  
\end{enumerate}
 
\ed

The following lemma is due to Plumstead, which is an adaptation of \cite[Example 4]{P}, tailored to our needs. This serves as one of the fundamental building blocks of this article. Here we provide a detailed proof. In the upcoming sections, some of the arguments used in the following proof are elided to avoid repeating the same arguments.

\bl\label{plgd}
Let $s$ be a non-zero divisor in $R$ such that $\dim(R_s)\le d-1$. Then there exists a generalized dimension function $d:\Spec(R)\to \mathbb{N}$ such that $d(\p)\le d-1$ for all $\p\in \Spec(R)$.

\el

\proof Let $\p\in \Spec(R)$. Then we note that either $s\in \p$ or $s\not \in \p$. That is, the prime ideal $\p\in  V(s)\cup \Spec(R_s)$, where $V(s)=\{\p\in \Spec(R):s\in \p\}$. Hence $\Spec(R)= S_1\cup S_2$, where $S_1=\Spec(R_s)$ and $S_2=V(s)$. We define two functions:
  \begin{enumerate}[\quad \quad (1)]
  	\item $d_1:S_1\to \mathbb{N}$ such that $d_1(\p)=\dim(R_s/\p)$ for all $\p\in S_1$;
  	\item  $d_2:S_2\to \mathbb{N}$ such that $d_2(\p)=\dim(R/\p)$ for all $\p\in S_2$.
  \end{enumerate}
  Let $I\subset R$ be an ideal. We observe that the equality $V(I)\cap S_1=V(I_s)$. Since $R$ is a Noetherian ring, it follows from the primary decomposition of $I_s$ that $V(I)\cap S_1$ has finitely many minimal elements with respect to the partial ordering $<<_1$ induced by $d_1$ as defined in Definition \ref{def}. On the other hand $V(I)\cap S_2=V(I)\cap V(s)=V(I+\langle s\rangle )$. Again, from the primary decomposition of $I + \langle s \rangle$, it follows that $V(I) \cap S_2$ also has finitely many minimal elements with respect to $<<_2$, induced by $d_2$. Therefore, both $d_1$ and $d_2$ are generalized dimension function on $S_1$ and $S_2$, respectively.
 
  Moreover, we notice that $d_1(\p)\le\dim(R_s)\le d-1$ for all $\p\in \Spec(R_s)=S_1$. On the other hand, since $s$ is a non-zero divisor, we obtain that $d_2(\p)\le \dim(R)-1\le d-1$ for all $\p\in S_2$.

Following \cite[Example 2]{P} we can define a generalized dimension function $d:\Spec(R)\to \mathbb{N}$ such that $d(\p)=d_1(\p)$ if $\p\in S_1$ and $d(\p)=d_2(\p)$ if $\p\in S_2$. Then we note that $d(\p)\le d-1$ for all $\p\in \Spec(R)$. This completes the proof.\qed

%\rmk Here we would like to point out an observation which was due to Plumstead. Suppose that $s\in R$ such that there exist generalized dimension functions $d_1:\Spec(R_s)\to \mathbb N$ such that $d_1(\p)\le d-1$ and  $d_2:\Spec(R/\langle s \rangle )\to \mathbb N$ such that $d_2(\p)\le d-1$. Then, one can obtain a generalized dimension function $d:\Spec(R)\to \mathbb N$ such that $d(\p)\le d-1$.
\smallskip

The next lemma has been used to reproduce Plumstead's results on patching in our setup. The following version is stated in a slightly more general form compare to the one in \cite[Lemma 2.17]{MPHIL}. The proof is essentially a verbatim copy of the same. Hence we skip the proof.

\bl \label{patchinglemma}
 Let $M$ be a $R[T]$-module. Let $N$ be a projective $R$-module. Let $\alpha(T)$ and $ \beta(T)$ be two surjective ($R[T]$-linear) maps from $M\surj N[T]$ such that $\alpha(0)=\beta(0)$, where $N[T]=N\otimes_R R[T]$. Furthermore, assume that the $R[T]$-modules $\ker(\alpha(T))$ and $\ker(\beta(T))$ are extended from the base ring $R$. Then there exists an automorphism $\Gamma(T)\in \Aut(M)$ such that (1) $\Gamma(0)=\text{Id}$ and (2) $\alpha(T)=\beta(T)\Gamma(T)$.
\el

The next theorem is an accumulation of several results, derived from a pivotal result originally due to Eisenbud and Evans \cite{EE}. This has been used extensively throughout the article. This version is recollected from \cite{P} and \cite[2.5]{SMBAR}.

\bt\label{eept}
  Let $\CP\subset\Spec(R)$ be a subset and $\delta:\CP \to \mathbb{N}$ a generalized dimension function. Let $M$ be a $R$-module such that $\mu_{\p}(M)\ge 1+\delta(\p)$ for all $\p\in \CP$. Let $(r,m)\in R\oplus M$ be a basic element on $ \CP$. Then there exists an element $m'\in M$ such that $m+rm'$ is a basic element on $\CP$. Consequently, if $M$ is a projective $R$-module such that $\rank(M_{\p})\ge d$ for all $\p\in \Spec(R)$ and $(r,m)\in \Um(R\oplus M)$, then there exists an element $m'\in M$ such that $\hh(\CO_M(m+rm'))\ge d$.

\et
The following lemma is due to Bhatwadekar and Roy \cite[Lemma 4.1]{SMBAR}. This has been crucially used in Proposition \ref{cansur}.

\bl\label{SMBAR}
Let $B\subset C$ be rings such that $\dim(B)= \dim(C)=d$. Assume that, there exists an element $x\in B$ such that $B_x=C_x$. Then 
\begin{enumerate}
	\item $B/\langle 1+xb\rangle B =C/\langle 1+xb\rangle C$ for all $b\in B$;
	\item if $\ma\subset C$ is an ideal such that $\hh(\ma)\ge d$ and $\ma+\langle x\rangle C=C$. Then there exists an element $b\in B$ such that $1+xb\in \ma$. 
\end{enumerate}
\el

We conclude this section with a lemma, which is analogous to \cite[Lemma 2]{P} in our setup. It enables us to patch two isomorphisms in the fiber product diagram used in this article. Before that, we define the following notation.

\smallskip

\notation 
Let $A$ be a subring of $R[T]$, that contains $R$. We define $\MI_A:=\{f(T)\in A:f(0)=0\}$. Then $\MI_A$ is an ideal in $A$. If $A$ is a geometric subring of $R[T]$ then one may observe that $\hh(\MI_A)\ge 1$.

\bl\label{pic}
Let $A$ be a geometric subring of $R[T]$. Let $s,t\in R$ such that $\langle s\rangle R+\langle t\rangle R=R$. Moreover, assume that  $A_s=R_s[T]$. Let $M$ and $M'$ be $A$-modules such that, there exist isomorphisms $\sigma_1:M_s\iso M_s'$ and $\sigma_2:M_t\iso M_t'$ with the following properties.

\begin{enumerate}[\quad \quad (1)]
	\item $(\sigma_1)_t\otimes (A/\MI_A)_{st}=(\sigma_2)_s\otimes (A/\MI_A)_{st}$;
	\item $M_{st}$ is a free module.
\end{enumerate} 
Then there exists an isomorphism $\sigma:M\iso M'$ such that $ \sigma\otimes (A/\MI_A)_s=(\sigma_1)\otimes (A/\MI_A)_{s}$ and $ \sigma\otimes (A/\MI_A)_t=(\sigma_2)\otimes (A/\MI_A)_{t}$.
\el
\proof We define $\phi(X):M_{st}[X]\iso M'_{st}[X]$ by $\phi(X):= (\sigma_1\otimes A_{st}[X])\circ [(\sigma_1)_{t}^{-1}\circ (\sigma_2)_s](XT)$, where $X$ is an indeterminate. Then $\phi(X)$ is an isomorphism satisfies the properties (i) $\phi(0)=(\sigma_1)_t$, (ii) $\phi(1)=(\sigma_2)_s$ and (iii) ${\phi(X)}\equiv \text{Id}\mod \MI_A R_{st}[T,X]$. Therefore, applying \cite[Lemma 1]{P} we get an isomorphism $\sigma:M\iso M'$ with the required property.\qed

\section{ Existence of basic elements}\label{3}

In this section, we prove an analogy of conjecture \ref{c1} regarding the existence of basic elements in certain modules over a geometric subring of a polynomial ring. We begin with the following lemma, which shows that it is sufficient to assume the module is torsion-free in order to prove our main theorem in this section. The proof is essentially contained in \cite[Theorem 2, paragraph 1]{P}. Hence we omit the proof.
\bl\label{tl}
 Let $M$ be a $R$-module and let {\setlength{\abovedisplayskip}{-1pt}
 	\setlength{\belowdisplayskip}{-1pt}$$M_{\text{tor}}=\{m\in M:\text{ there exists a non-zero divisor } s\in R \text{ such that } sm=0\}$$} be the torsion submodule of $M$. Let $M'=M/M_{\text{tor}}.$ If $M'$ has a basic element, then $M$ has a basic element.
\el
The next proposition is an improvement of \cite[Proposition 1]{P} in our setup.

\bp\label{pobe}
Let $A$ be a geometric subring of $R[T]$. Let $s,t\in R$ such that $A_s=R_s[T]$ and $\langle s\rangle R+\langle t\rangle R=R$.  Let $M$ be an $A$-module. Let ``bar'' denote going modulo $\MI_AM$. Let $z$ be an element in $\ol M$. Suppose that, there exist basic elements $m_1\in M_s$ and $m_2\in M_t$ such that $\ol m_1=(z)_s$ and $\ol m_2=(z)_t$. Let $$N_1=\frac{M_{s}}{m_1A_{s}} \text{ and }N_2=\frac{M_t}{m_2A_t}.$$ Furthermore, assume that $(N_1)_t$ and $(N_2)_s$ are extended projective modules from the base ring $R_{st}$. Then there exists an basic element $m\in M$ such that $\ol m=z$.

\ep
\proof  Since $(N_1)_{t}$ and $(N_2)_s$ are projective modules over $A_{st}$, we have the following identity. 
\begin{equation}\label{equ:A}\tag{A}
	M_{st}\cong m_1A_{st}\oplus (N_1)_{t}\cong m_2A_{st}\oplus (N_2)_s
\end{equation}
First we observe that $\ol{( m_1)_{t}}=\ol{( m_2)_{s}}=(z)_{st}$. Therefore, the definition of $N_i$ ($i=1,2$) give us the equality $\ol{( N_1)_{t}}=\ol{( N_2)_s}$, which we call $N$. Since both the modules $( N_1)_{t}$ and $(N_2)_s$ are extended from the base ring $R_{st}$, there exist isomorphisms $\phi_1(T):( N_1)_{t}\iso N[T]$ and $\phi_2(T):( N_2)_s\iso N[T]$ such that $\phi_1(0)=\phi_2(0)=\text{Id}$. We define the following surjections.
\begin{enumerate}[\quad \quad (a)]
	\item $\eta_1:A_s\surj m_1A_s$ such that $\eta_1(1)=m_1$;
	\item  $\eta_2: A_t\surj m_2A_t$ such that $\eta_2(1)=m_2$.
\end{enumerate}
From \ref{equ:A}, it follows that the maps $(\eta_1)_t$ and $(\eta_2)_s$ are in fact isomorphisms. Moreover, $(\eta_1)_t$ and $(\eta_2)_s$ will canonically induce maps, say $\widetilde{\eta_1}:A_{s}\to M_{s}$ and $\widetilde{\eta_2}:A_{t}\to M_{t}$ such that $(\widetilde{\eta_1})_t$ and $(\widetilde{\eta_2})_s$ are injective maps. Now we consider the following commutative diagram of short exact sequences
\begin{equation}\label{equ:B}\tag{B}
	\begin{tikzcd}
		0\ar[r] & R_{st}[T]\ar[rr,"(\widetilde{\eta_1})_t"]\arrow[equal]{d} && M_{st}\ar[rr,"\rho_1(T)"]\arrow[d,"\Gamma(T)",dashed] && N[T]\ar[r] \ar[d,equal] & 0\\
		0\ar[r] & R_{st}[T]\ar[rr,"(\widetilde{\eta_2})_s"]&& M_{st}\ar[rr,"\rho_2(T)"] && N[T]\ar[r] & 0
	\end{tikzcd} 
\end{equation}
where $\rho_i(T)$, $i=1,2$, are defined as the compositions of the canonical surjections induced from the identity (\ref{equ:A}), along with the isomorphisms $\phi_i(T)$ $(i=1,2)$. Moreover, from the definition of $\rho_i(T)$, it follows that $\rho_1(0)=\rho_2(0)$. Therefore, using Lemma \ref{patchinglemma}, there exists an automorphism $\Gamma(T)\in \Aut(M_{st})$ with $\Gamma(0)=$ Id such that $\rho_1(T)=\rho_2(T)\circ \Gamma(T)$. Using Quillen's splitting lemma \cite[Theorem 1]{Q}, the automorphism $\Gamma(T)$ splits. We consider the following fiber product diagram.
$$\begin{tikzcd}
	A \ar[rrr] \ar[dd]\ar[dr, dashed, "\eta"] &&&  A_s \ar[dd]\ar[dr,"\widetilde{\eta_1}" ] \\
	& M \ar[rrr,crossing over, dashed]  &&& M_s\ar[dd] \\
	A_{t} \ar[rr] \ar[dr,"\widetilde{\eta_2}"] & &   A_{st} \ar[dr] \ar[r,equal] & A_{st}\ar[dr] \\
	& M_{t}\ar[rr]  \arrow[uu,crossing over, leftarrow,dashed] && M_{st} \arrow{r}{\Gamma^{-1}}[swap]{\sim} & M_{st}
\end{tikzcd}$$ 
Since $\Gamma(T)$ splits, from the universal property of the fiber product we obtain a unique map $\eta:A\to M$ such that $(\eta)_s=\widetilde{\eta_1}$ and $(\eta)_t=\widetilde{\eta_2}$. Let $m=\eta(1)$. Since $\Spec(A)=\Spec(A_s)\cup \Spec(A_t)$, a local checking ensures that $m$ is a basic element of $M$. Furthermore, we have $\ol m=z$, as again it holds locally. This completes the proof. \qed

We are now ready to give the proof of the main theorem in this section.

\bt\label{ebe}
Let $A$ be a geometric subring of $R[T]$. Let $M$ be an $A$-module such that 
$	\mu_{\p}(M)\ge \dim(A/\p)$ for all minimal prime ideal $\p\in \Spec(A)$. Then $M$ has a basic element. 
\et

\proof

Let $\eta\subset A$ be the nilradical of $A$. Since $\Spec(A/\eta)=\Spec(A)$, every basic element of $M/\eta M$ lifts to a basic element of $M$. Therefore, without loss of generality, we may assume that $A$ is reduced. In particular, this implies $R$ is reduced. Moreover, it follows from Lemma \ref{tl} that it is enough to assume $M$ is a torsion-free module. Since $A$ is a geometric subring of $R[T]$, there exists a non-zero divisor $s_1\in R$ and $k\ge 0$ such that $s_1^kT\in A$. We take $f(T)=s_1^kT$, then $f\in \MI_A$ is a non-zero divisor in $A$. Throughout the proof, the notation ``bar'' denotes going modulo $\langle f \rangle M$. Applying Theorem \ref{eept} on the module $\ol{M}$, we obtain a basic element $z\in \ol M$. The remaining part of the proof is given in steps.\\
\textbf{Step - 1.} In this step, we show the existence of a basic element in $M_s$ for some suitably chosen $s\in R$. Let $\MT$ be the multiplicatively closed set of all non-zero divisors in $R$. Then we notice that $\MT^{-1}A=\MT^{-1}R[T]=(k_1\times ...\times k_n)[T]$. Since $M$ is a torsion-free module, the module $\MT^{-1}M$ is a free module over $\MT^{-1}A$. Also one may observe that, the canonical image of $z$ in $\ol{\MT^{-1}M}$, say $(z)_{\MT} ,$ is a basic element in $\ol{\MT^{-1}M}$. Let $m_1'$ be a lift of $(z)_{\MT}$ in $\MT^{-1}M$. Due to the fact that $\dim(\MT^{-1}A)=1$, one may apply Theorem \ref{eept} on $(m_1',T)$, to obtain a basic element $m_1:=m_1'+Tm_1''$ in $\MT^{-1}M$, for some suitably chosen $m_1''\in \MT^{-1}M$. Here an observation is that, $\ol m_1=(z)_{\MT}$ and the module $\MT^{-1}M/m_1\MT^{-1}A$ is a free $\MT^{-1}A$-module. Since $M$ is a finitely generated module, there exists a non-zero divisor $s_2\in R$ such that the modules $M_{s_2}$ and $M_{s_2}/m_1A_{s_2}$ are free $A_{s_2}$-modules and $m_1\in \Um (M_{s_2})$. Let $s=s_1s_2$. Then in the ring $A_s$ we have the following.
\begin{enumerate}[\quad \quad (1)]
	\item $s\in R$ is a non-zero divisor;
	\item $A_s=R_s[T]$;
	\item $m_1\in \Um(M_s)$ such that $\ol m_1=(z)_s$;
	\item both the modules $M_s$ and $M_s/m_1A_s$ are free $A_s$-modules.
\end{enumerate}	
\textbf{Step - 2.} Let  $S=\{1+sr:r\in R \}$. This step is devoted to establish the existence of a basic element in $S^{-1}M$. Let us denote $B=S^{-1}A$ and $L=S^{-1}M$. Then we note that, $B$ is a geometric subring of $(S^{-1}R)[T]$, which we write $R_{1+sR}[T]$. Since $s$ is a non-zero divisor, using Lemma \ref{plgd} we get a generalized dimension function $d:\Spec(B)\to \mathbb{N}$ such that $d(\p)\le d$ for all $\p\in \Spec(B)$. We choose a lift, say $m_2'\in L$, of $(z)_S$. One may observe that, for any $\p\in \Spec(B)$, if $f \in \p$, then $m'_2\not \in \p L_{\p}$, as $\ol{\p}\in \Spec(\ol B)$. This eventually implies that, the element $(f,m_2')$ is a basic element of $B\oplus L$. Applying Theorem \ref{eept} on the element $(f,m_2')$, we obtain a basic element $m_2:=m_2'+fm_2''$, of $ L$, where $m_2''\in L$.\\ 
\textbf{Step - 3.}  We devote the rest of the proof in preparation to apply Proposition \ref{pobe}, which allows us to patch $m_1$ and $m_2$ together to obtain a basic element $m\in M$. We define the following modules.
$$N_1:=\frac{M_{s}}{m_1A_{s}} \text{ and }N'_2:=\frac{L}{m_2B}.$$  
From (4), it follows that the module $N_1$ is an free $A_s$-module.  Since $M_{s}$ is a free $A_s$-module and $(m_2)_s\in \Um(L_s)$, the module $(N'_2)_s$ is a stably free module over the ring $R_{s(1+sR)}[T]$. Let $\mq\in \Spec(B_{s})=\Spec(R_{s(1+sR)}[T])$ be a minimal prime ideal. Then we note that, $A\cap \mq\in \Spec(A)$, is also a minimal prime ideal in $A$. Moreover, as $(\mq\cap A) \cap S=\phi$, we obtain the equality  $(L_{s})_{\mq}=M_{A\cap \mq}$ of modules. Therefore, $$\mu_{\mq}({(N'_2)_s})=\mu_{\mq}(L_{s})-1=\mu_{A\cap \mq}(M)-1\ge \dim(A/A\cap \mq)-1.$$ 
One may observe that $ \dim(A)-1\ge \dim(B_{s})$. 
Implies that $\dim(A/A\cap \mq)-1\ge \dim(B_{s}/ \mq)$. Applying cancellation result, as stated in \cite[Corollary 1]{P}, on the ring $R_{s(1+sR)}[T]$, we obtain that the module $(N'_2)_s$ is in fact a free module. 

Since all modules considered in the theorem are finitely generated, there exists $t\in S$ such that $m_2$ is a basic element of $M_t$ such that $\ol m_2=(z)_t$ and if we take $N_2=\frac{M_t}{m_2A_t}$, then $(N_2)_s$ is an extended projective module from the base ring $R_{st}$. Now applying Proposition \ref{pobe}, we obtain a basic element $m\in M$ such that $\ol m=z$. \qed
\smallskip

	\section{Cancellation of projective modules}\label{4}
	
	In this section, we study the cancellation problem for projective modules over a geometric subring of $R[T]$. We begin with the following technical proposition, similar to \cite[Remark 2.5]{m82}. This has been used frequently throughout this section to establish various cancellation results.
	  \smallskip
	  
	  \notation Let $M$ be an $R$-module. Given an automorphism $\alpha\in \Aut(R\oplus M)$ and an element $(f,m)\in R\oplus M$, we denote $(f,m)\alpha$ as the image of $(f,m)$ under $\alpha$. This notation is introduced for consistency with the case when $M$ is a free module over $R$.

	\bp\label{canl}
	Let $A$ be a reduced geometric subring of $R[T]$. Let $M$ be a torsion-free $A$-module such that $M_g$ is a projective $A_g$-module of rank $\ge d+1$, for some non-zero divisor $g\in \MI_A$. Let $(f,m)\in \Um(A\oplus M)$ such that
	\begin{enumerate}[\quad \quad (1)]
		\item $f-1\in \langle g\rangle A$ and $m\in \langle g\rangle M$; \label{14.2}
		\item $g+g^2g'\in \CO_M(m)$, for some $g'\in A$.\label{24.2}
	\end{enumerate}
	Then $\frac{A\oplus M}{(f,m)A}\cong M$. Moreover, there exists an $\alpha\in \Aut(A\oplus M)$ such that $(f,m)\alpha=(1,0)$. Furthermore, if $M$ is a projective module, then the conclusion holds even without (\ref{24.2}).
	\ep
	\proof  First we notice that if we identify $\frac{A\oplus M}{(1,0)A}$ simply with $M$, then we are not losing any generality. Some more observations are as follows: suppose that there exists an $\alpha\in \Aut(A\oplus M)$ such that $(f,m)\alpha=(1,0)$. Then we obtain the required isomorphism induced by $\alpha$ from the following short exact sequence
	\begin{equation*}
		\begin{tikzcd}
			0 \arrow{r} & A \arrow[rr,"{1\to (f,m)}"] && A\oplus M\arrow["\alpha'"]{rr}  && \frac{A\oplus M}{(1,0)A}=M\arrow{r} & 0,
		\end{tikzcd}
	\end{equation*}
	where $\alpha'$ is the composition of $\alpha$ along with the canonical map $A\oplus M\surj \frac{A\oplus M}{(1,0)A}$. From (\ref{14.2}) there exists $h\in A$ such that $f=1+gh$. Let $\gamma:M\to \frac{ A\oplus M}{(f, m)A}$ be the canonical map. Let $M'=\frac{ A\oplus  M}{( f, m)A}$ and let ``bar'' denote going modulo $\MI_A$ as well as $\MI_AM$. Then we notice that $\ol \gamma=\text{Id}$. The remaining part of the proof is divided into two steps.\\
	% \begin{equation}\tag{B}
		%\gamma:\frac{\ol A\oplus \ol M}{(\ol f,\ol m)A}\xrightarrow{\sim}  \ol M.
		% \end{equation} 
	\textbf{Step - 1.} 
	 We start with the observation that, since $(f,m)\in \Um(A\oplus M)$, we have $A\oplus M'\cong A \oplus M$. Therefore, the module $M'$ is torsion-free, as it is a direct summand of a torsion-free module $A\oplus M$. Following the arguments given in [Theorem \ref{ebe}, Step - 1, paragraph 1] we can find a non-zero divisor $s\in R$ such that 
	\begin{enumerate}[\quad \quad (a)]
		\item $A_s=R_s[T]$;\label{a4.2}
		\item $M_s$ and $M_s'$ are free modules over the ring $A_s$.\label{b4.2}
		%	\item $\sigma_1\in \E(A_s\oplus M_s)$ such that $(f,m)\sigma_1=(1,0)$.
	\end{enumerate} 
	In particular, from (\ref{b4.2}) it follows that, the modules $M_s$ and $M_s'$ are extended from the base ring $R_s$. Moreover, since $ \ol M_s=\ol M_s'$, there exists an isomorphism $\sigma_1:M_s\iso  M_s'$ such that $\ol{\sigma_1}=\text{Id}$.\\
	\textbf{Step - 2.} 
	Let $S=\{1+sr:r\in R\}$. We denote $B=S^{-1}A$, $L=S^{-1}M$ and $L'=S^{-1}M'$. Since $s$ is a non-zero divisor, using Lemma \ref{plgd} we can obtain a generalized dimension function $d:\Spec(B)\to \mathbb{N}$ such that $d(\p)\le d$ for all $ \p\in \Spec(B)$.  Let $\CS= \{\p\in \Spec(B):g\not \in \p\}$ and $\delta=d|_{\CS}$, be the restriction of $d$ on $\CS$. Then $\delta:\CS\to \mathbb{N}$ is a generalized dimension function [cf. the proof of Lemma \ref{plgd}, paragraph 3] such that  $\delta(\p)\le d$ for all $ \p\in \CS$. 
	
	We claim that, the element $(g^2f,m)$ is a basic element in $B\oplus L$ on $\CS$. To establish our claim, we note that it is enough to show that $\CO_{B\oplus L}(g^2f,m)\not\subset \p$, for any $\p\in \CS$ which contains $f$. Let $\p\in \CS$ be such that $f\in \p$. Since $(f,m)\in \Um(B\oplus L)$, there exists $\widetilde{\lambda}\in (B\oplus L)^*$ such that $\widetilde\lambda(0,m)\not\in \p$. In particular, this implies $\CO_{B\oplus L}(g^2f,m)\not\subset \p$.
	
	Applying Theorem \ref{eept} on the triplet $\{(g^2f,m),\CS, B\oplus L\}$, we get a basic element $m''=m+m'fg^2\in L$ on $ \CS$, where $m'\in L$. From the definition of $m''$ and using hypothesis (\ref{14.2}), it follows that $m''\in S^{-1}(\langle g\rangle M)$. As $M_g$ is a projective $A_g$-module, the element $m''\in \Um(L_g)$. Therefore, there exists an integer $k\ge 1$ such that $g^k\in \CO_{L}(m'')$. Using hypothesis (\ref{24.2}) we obtain that $g^{k-1}\in \CO_{L}(m'')$. Repeating this argument finitely many times we will eventually get $g\in \CO_{L}(m'')$.
	
	Let $\lambda\in L^*$ such that $\lambda(m'')=g$. In the remaining part of this step we construct a series of automorphisms of $B\oplus L$, which eventually transform $(f,m)$ to $(1,0)$. For each $i=1,2,3$ and $(F,p)\in (B\oplus L)$ we consider the following automorphisms of $B\oplus L$.
	\begin{enumerate}[\quad \quad $\bullet$ ]
		\item $(F,p)\phi_1=(F, p+Fg^2m')$
		\item $(F,p)\phi_2= (F-h\lambda(p),p)$
		\item $(F,p)\phi_3=(F,p-Fm'')$
	\end{enumerate} 
	We define $\chi:=\phi_3\circ \phi_2\circ\phi_1$. It follows from the explicit computation of $\chi$ that $(f,m)\chi=(1,0)$. Since $m''\in S^{-1}(gM)$ for all $(F,p)\in B\oplus L$, one may observe that 
	\begin{equation}\tag{A}\label{B1}
		(\ol F,\ol p)\ol\chi= (*,\ol p).
	\end{equation} 
	We consider the following commutative diagram
	$$
	\begin{tikzcd}
		0\ar[r] & B\ar[rr,"{1\to (1,0)}"]\arrow[equal]{d} && B \oplus L\ar[rr]\arrow[d,"\chi^{-1}"] && L\ar[r] \ar[d,"{\sigma_2'}",dashed] & 0\\
		0\ar[r] & B\ar[rr,"{1\to (f,m)}"] && B\oplus L\ar[rr] && L'\ar[r] & 0
	\end{tikzcd} 
	$$
	where $\sigma_2'$ is the canonical isomorphism induced by $\chi^{-1}$, which makes the above diagram commutative. Furthermore, from (\ref{B1}) it follows that $\ol \sigma_2'=\text{Id}.$ Since all modules considered in this proof are finitely generated, there exists $t\in S$ such that we can find an isomorphism $\sigma_2:M_t\iso M_t'$ with the property that $\ol{\sigma}_2=\text{Id}$. Therefore, applying Lemma \ref{pic} we get the required isomorphism $\sigma:M\iso M'$ such that $\ol \sigma=\text{Id}.$ 
	
We now consider the following commutative diagram
	\begin{equation}\tag{B}\label{C1}
		\begin{tikzcd}
			0\ar[r] & A\ar[rr,"{1\to (f,m)}"]\arrow[d,equal] && A \oplus M\ar[rr,"\pi_1"]\arrow[d,"\alpha", dashed] &&  \frac{ A\oplus M}{(f, m)A}=M'\ar[r] \ar[d,"{\sigma}^{-1}"] & 0\\	0\ar[r] & A\ar[rr,"{1\to (1,0)}"] && A\oplus M\ar[rr,"\pi_2"] &&  \frac{ A\oplus M}{(1, 0)A}=M\ar[r] & 0
		\end{tikzcd} 
	\end{equation}
	where $\pi_i$ are the canonical quotient maps, $i=1,2$. Since the exact sequences split, there exists an $\alpha\in \Aut(A\oplus M) $ such that $(f,m)\alpha=(1,0)$. This concludes the proof. \qed

	We are now ready to prove the cancellation result for projective modules over a geometric subring of $R[T]$. The idea of the proof is to establish the fact that, in Proposition \ref{ebe}, if $M$ is a projective module, then even hypothesis (\ref{14.2}) is not required. 
	
%	Here we would like to mention that, the following cancellation result is proved in \cite{RSP} for Rees algebras. However, I could not understand their arguments, especially \cite[Lemma 3.7]{RSP}.
	
	\bt\label{can}
	Let $A$ be a geometric subring of $R[T]$. Let $P$ be a projective $A$-module such that $\rank(P)\ge d+1$. Then the projective module $P$ is cancellative.
	\et
	
	\proof We first note some general observations: let $\eta$ be the nilradical of $R$ and let $R_{\text{red}}=R/\eta$ and $A_{\text{red}}=A/\eta A$. Then the ring $ A_{\text{red}}$ is also a geometric subring of $R_{\text{red}}[T]$. Moreover, since we can lift automorphisms up to nilradical, we may assume without loss of generality that $R$ is reduced. Let $\psi:A\oplus Q\iso A\oplus P$ be an isomorphism such that $(1,0)\psi=(f,m)$. Let $s$ be a non-zero divisor in $R$ such that $A_s=R_s[T]$. Then to establish the theorem, it is enough to find $\alpha\in \Aut(A\oplus P)$ such that $(f,m)\alpha=(1,0)$. The remaining part of the proof is devoted in showing this.
	
	We note that $(f,m)\in \Um(A\oplus P)$. From our choice of the ring $A$, there exists $n\in \mathbb{N}$ such that $s^nT\in A$. Let $g(T)=s^nT$. Then apparently $g\in \MI_A$ is a non-zero divisor. Let ``tilde" denote going modulo $\langle g\rangle A$ as well as $\langle g\rangle P$. Since $g$ is a non-zero divisor, using Theorem \ref{eept}, there exists a transvection $\widetilde{\tau} \in \E(\widetilde{A}\oplus\widetilde{P})$ such that $(\widetilde{f},\widetilde{p} )\widetilde{\tau}=(\widetilde{1},\widetilde{0})$. Recall that, we can always lift transvections of $\widetilde{A}\oplus\widetilde{P}$ to an automorphism of ${A}\oplus {P}$. Let $\tau\in \Aut(A\oplus P)$ be a lift of $\widetilde{\tau}$. Therefore, by altering $(f,m)$ with $(f,m)\tau$, we may further assume that $m\in \langle g\rangle M$ and $f=1+gh$ for some $h\in A$. Hence we may apply Proposition \ref{canl} to conclude the proof.\qed

In the remaining part of this section, we aim to extend Theorem \ref{can} to a geometric subring of $R[T,f^n]$, where $f\in R[T]\setminus R$ is a non-zero divisor and $n\in \mathbb Z$. To achieve this, we rely on the following proposition, which allows us to reduce the problem to a problem over a geometric subring of $R[T]$. The proof of the proposition draws inspiration from \cite[Theorem 4.3]{SMBAR}.
	
	\bp \label{cansur}
	Let $B\subset C$ be reduced rings of dimension $d\ge 1$ and let $x\in B$ be a non-zero divisor in $C$ such that $B_x=C_x$. Let $P$ be a projective $C$-module of rank $\ge d$ and $(c,p)\in \Um(C \oplus P)$. Then there exist (a) a transvection $\tau\in \E(C \oplus  P)$, (b) a torsion-free $B$-module $M$ and (c) a unimodular element $(b,m)\in \Um(B\oplus M)$ such that:
		\begin{enumerate}
			\item $M_x\cong P_x$;
			\item $(c,p)\tau =(b,m)$;
			\item $b-1\in \langle x\rangle B$ and $m\in \langle x\rangle M$;
			\item $x+x^2g\in \CO_{ M}(m)$, for some $g\in B$.
		\end{enumerate}

	\ep
	\proof Since $x$ is a non-zero divisor in $C$ and $\rank(P)\ge d$, it follows from \cite{S} that $P/\langle x\rangle P$ has a unimodular element. Therefore, by altering $(c,p)$ up to a transvection of $C\oplus P$, we may assume without loss of generality that $c\in \langle x \rangle C$. Moreover, applying Theorem \ref{eept}, we may replace $p$ with $p+cq$ for some $q\in P$ and further assume that $\hh(\CO_{ P}(p))\ge d$. Let $\ma=\CO_{ P}(p)$. Because of $(c,p)\in \Um(C\oplus P)$, we get $\langle x \rangle C+\ma=C$. Therefore, using Lemma \ref{SMBAR}, we obtain an element $1+xy\in \ma$, for some $y\in B$. Hence one can choose an $\alpha\in P^*$ such that $\alpha(p)=1+xy$. Let $c=xh'$ for some $h'\in C$. Since $\frac{B}{\langle 1+xy\rangle B}=\frac{C}{\langle 1+xy\rangle C}$, there exists $b'\in B$ such that $h'=b'+\lambda(1+xy)$, for some $\lambda\in C$. Then we have 
	\begin{equation}\tag{i}\label{i}
	c+[1+x(1-\lambda-y-b')](1+xy)=(1+x+b''x^2),
\end{equation}
where $b''=y-yb'-y^2\in B$. We take $b=1+x+b''x^2$. Then one may notice that $b\in B$. From (\ref{i}) it follows that, the unimodular element $(c,p)$ can be translated to $(b,p)$ by applying a transvection of $C\oplus P$. Moreover, if we take  $m=p-bp$, then $(b,m)\in \Um(C\oplus P)$ and $(b,p)$ can be translated to $(b,m)$ by applying a transvection of $C\oplus P$. Furthermore, as a consequence we have 
	\begin{equation}\tag{ii}\label{ii}
		b-1\in \langle x \rangle B, \text{ and } m\in \langle x\rangle P.
\end{equation}
Let $m=xp_1$ for some $p_1\in P$. Then we note that \begin{equation}\tag{iii}\label{iii}
	\alpha(m)=(-x-b''x^2)(1+xy)=-x+x^2h,
\end{equation} for some $h\in B$. This implies $\alpha(p_1)=-1+xh\in B$. Let $\{p_1,p_2',...,p_k'\}$ be a set of generators of $P$ and $\alpha(p_i')=c_i\in C$ for all $i=2,...,k$. Using the identity $	\frac{B}{\langle -1+xh \rangle B}=\frac{C}{\langle -1+xh\rangle C} $ we can find $\lambda_i\in C$ and $b_i\in B$ such that $c_i=b_i+\lambda_i(-1+xh)$ for all $i=2,...,k$. We define $p_i:=p_i'-\lambda_ip_1$, where $i=2,...,k$. Then one can notice that $\{p_1,p_2,...,p_k\}$ generates $P$ as a $C$-module such that $\alpha(p_i)\in B$ for all $i=1,...,k$.
	
	We set $M=\sum_{i=1}^{k}Bp_i$. Then we claim that $M$ is a torsion-free module over $B$. To see this, first we observe that, the element $x$ is also a non-zero divisor in $B$. Let $m\in M\setminus \{0\}$ such that $am=0$, where $a\in B\setminus \{0\}$. Since $M$ is a $B$-submodule of a torsion-free module $P$, the element $a$ must be a zero-divisor in $C$. Then there exists a minimal prime ideal $\p\in \Spec(C)$ such that $a\in \p$. Since $B_x=C_x$ and $x$ is a non-zero divisor in both $B$ and $C$, there exists a canonical one-to-one correspondence between the minimal prime ideals of $B$ and $C$. In particular, this implies $a$ is a zero-divisor in $B$. This proves our claim. 
	
	It follows from the construction of $M$ that the canonical map, say $\Gamma:M\otimes_B C\surj P$, is a surjection. Moreover, if we take $\omega\in \ker(\Gamma)$, then there exists $k\ge 0$ such that $x^k\omega=0$. Therefore, the map $\Gamma_x$ is an isomorphism. This implies $M_x\cong P_x$. One may treat $\alpha$ as an element of $M^*$ by taking restriction on $M$. Since $(b,m)\in \Um(C\oplus P)$, it follows from our choice of $M$ that $(b,m)\in \Um(B_x\oplus M_x)$. Therefore, for some integer $t\ge 0$, the element $x^{t}\in \CO_{B\oplus  M}(b,m)$. As some suitable power of $b \,a(=1+x+b''x^2)$, is co-maximal with $x^{t}$, the element \begin{equation}\tag{iv}\label{iv} (b,m)\in \Um(B\oplus M).
	\end{equation} Combining (\ref{ii}), (\ref{iii}) and (\ref{iv}), it follows that $(b,m)$ satisfies all the desired conditions.\qed

	The next lemma is an interesting consequence of the arguments used in Proposition \ref{cansur}. Before that, we define the following notation.
	\smallskip
	
\notation
 Let $u,v\in \Um_n(R)$. We define $ u \sim_{\E_n(R)} v$ if there exists an $\epsilon\in \E_n(R)$ such that $u\epsilon=v$, where  $\E_n(R)$ is the subgroup of $\GL_n(R)$, generated by all the elementary matrices.

	\bl\label{ic}
	Let $B$, $C$ and $x$ be as in Proposition \ref{cansur}. Then the canonical map {\setlength{\abovedisplayskip}{5pt}
		\setlength{\belowdisplayskip}{2pt}$$\begin{tikzcd}
		\frac{\Um_{d+1}(B)}{\E_{d+1}(B)}\arrow[r,  twoheadrightarrow] & \frac{\Um_{d+1}(C)}{\E_{d+1}(C)}
	\end{tikzcd}$$} is a surjection.
	\el
	\proof Let $(c,p)\in \Um_{d+1}(C)$, where $p=(c_1,...,c_d)\in C^d$. We follow the proof of Proposition \ref{cansur} till we obtain that $(c,p)\sim_{\E_{d+1}(C)}(b,p)$ with $b=1+x+b''x^2\in B$, where all the notations defined here are consistence with the notation used in the proof of Proposition \ref{cansur}, till that point. To establish the lemma we notice that it is enough to find some $(b_1,...,b_d)\in B^d$ such that (1) $(b,p)\sim_{\E_{d+1}(C)}(b,b_1,...,b_d)$ and (2) $(b,b_1,...,b_d)\in \Um_{d+1}(B)$. Let $z=1+b''x\in B$. Then note that $b=1+xz$. Using the identity $	\frac{B}{\langle 1+xz \rangle B}=\frac{C}{\langle 1+xz\rangle C} $ we can find $\mu_i\in C$ ( $i=1,...,d$ ) such that $c_i-\mu_i(1+xz)\in B$. Now if we take $b_i=c_i-\mu_i(1+xz)$ for all $i=1,...,d$, then it follows that \begin{equation}\tag{A}\label{A2}
		(b,p)	\sim_{\E_{d+1}(C)}(b,b_1,...,b_d).
	\end{equation} It is only remains to show that $(b,b_1,...,b_d)\in \Um_{d+1}(B)$. We observe that, from (\ref{A2}) it follows that $(b,b_1,...,b_d)\in \Um_{d+1}(C_x)=\Um_{d+1}(B_x)$. That is, for some $t\ge 0$, we must have $x^t\in \langle b,b_1,...,b_d \rangle B$. Now from our choice of $b$ it follows that $1\in \langle b,b_1,...,b_d \rangle B$. This completes the proof of the lemma.\qed 

We end this section with the following theorem, which is essentially an extension of Theorem \ref{can} to a geometric subring of $R[T,f^n]$, where $f\in R[T]\setminus R$ is a nonzero divisor and $n\in \mathbb{Z}$.

		\bt\label{cane}
		Let $f\in R[T]\setminus R$ be a non-zero divisor. Let $A$ be a geometric subring of $R[T,f^n]$, where $n\in \mathbb Z$. Let $P$ be a projective $A$-module such that $\rank(P)\ge d+1$. Then the projective module $P$ is cancellative.
		\et
		\proof Without loss of generality we may assume that $R$ is reduced. Let $\psi:A\oplus Q\iso A\oplus P$ be an isomorphism such that $(1,0)\psi=(c,p)$. Then $(c,p)\in \Um(A\oplus P)$. For all $n \geq 0$ the result follows from Theorem \ref{can}, so we assume that $n < 0$. Let $B=A\cap R[T]$. Then $B$ is a geometric subring of $R[T]$. Furthermore, we observe that there exists an element $x\in \MI_B$, which is a non-zero divisor in $A$ such that $A_x=B_x=R_s[T,\frac{1}{Tf}]$. Therefore, applying Proposition \ref{cansur} on $(c,p)$ (taking $C=A$) one may obtain (i) a torsion-free $B$-module $M$, (ii) a transvection $\tau\in \E(C \oplus  P)$ and (iii) $(f,m)\in \Um(B\oplus M)$ such that:
		\begin{enumerate}[\quad \quad (1)]
			\item $M_x\cong  P_x$;
	\item $(c,p)\tau=(f,m)$; 		\label{2c}
			\item $f-1\in \langle x\rangle B$ and $m\in xM$;
			\item $x+x^2h\in \CO_{ M}(m)$, for some $h\in B$.
		\end{enumerate}
		Here we note that, since the kernel of the canonical surjection $A\oplus (M\otimes_B A)\surj A\oplus P$ is the $x$-torsion submodule of $A\oplus(M\otimes_B A),$ any automorphism $\omega$ of $B\oplus M$ such that $(f,m)\omega=(1,0)$ will canonically induce an automorphism of $A\oplus P$ which takes $(f,m)$ to $(1,0)$. Therefore, in view of (\ref{2c}), finding such an $\omega\in \Aut(B\oplus M)$ will conclude the proof. Since $x$ is a non-zero divisor and $\rank(M_x)=\rank(P_x)\ge d+1$, one may apply Proposition \ref{canl} on $(f,m)$ to obtain the required automorphism. This concludes the proof. \qed

		\section{Efficient generation of modules}\label{5}
		In this section, we study the upper bound on the number of generators of a module, suggested by Eisenbud-Evans, over a geometric subring of $R[T]$. We begin with a lemma which is due to Mohan Kumar. The proof is essentially contained in \cite[Theorem 1, Step - 1]{NMK}. Hence we skip the proof to avoid repeating similar arguments. 
		\bl\label{l1}
		Let $M$ be a $R$-module. Let $M'=M/M_{tor}$. Suppose that $M'$ is generated by $e(M)$ many elements. Then $M$ is also generated by $e(M)$ many elements.
		\el
\notation With the same notations as of Lemma \ref{l1}, we denote $f(M)=\sup\{\mu_{\p}(M)+\dim(R/\p):\p\in \Spec(R)\text{ such that } M_{\p}\not=0\}$. Then it has been shown in \cite[Theorem 4.1, paragraph 2]{M} that $e(M)\ge f(M)-1$.

\smallskip
		
The next proposition is a refinement of \cite[Proposition 2]{P} in our setup.
		
\bp\label{pog}
	Let $A$ be a geometric subring of $R[T]$. Let $s,t\in R$ such that $A_s=R_s[T]$ and $\langle s\rangle R+\langle t\rangle R=R$.  Let $M$ be an $A$-module such that $M_{st}$ is an extended projective module from the base ring $R_{st}$. Let {\setlength{\abovedisplayskip}{5pt}
		\setlength{\belowdisplayskip}{5pt}$$\phi_1:(A_s)^m\surj M_s \text{ and } \phi_2:(A_t)^m\surj M_t$$} be two surjections such that 
	{\setlength{\abovedisplayskip}{5pt}
		\setlength{\belowdisplayskip}{5pt}$$(\phi_1)_t\otimes (A/\MI_A)_t=(\phi_2)_s\otimes (A/\MI_A)_s.$$ }Suppose that, the modules $(\ker(\phi_1))_{t}$ and $(\ker(\phi_2))_s$ are extended from the base ring $R_{st}$. Then there exists a surjection $\phi:A^m\surj M$ such that $\phi_s=\phi_1$ and $\phi_t=\phi_2.$
	\ep	
	\proof 	Let $\ker(\phi_i)=L_i,$ for $i=1,2$ and ``bar'' denote going modulo $\MI_AM$. We take $N=\ol M$. Since $M_{st}$ is an extended module from the base ring $R_{st}$, there exists an isomorphism $\Gamma(T):M_{st}\iso N_{st}[T]$ such that $\Gamma(0)=\text{Id}$. Now we have the following commutative diagram
	\begin{equation}\tag{A}\label{A3}
		\begin{tikzcd}
			0\ar[r] & (L_1)_{t}\ar[rr]\arrow{d} && (R_{st}[T])^m\arrow[rr,"\psi_1"]\arrow[d,dashed,"\omega(T)"] && N_{st}[T]\ar[r] \ar[equal]{d} & 0\\
			0\ar[r] & (L_{2})_s\ar[rr]&& (R_{st}[T])^m\arrow[rr,"\psi_2"]  && N_{st}[T]\ar[r] & 0
		\end{tikzcd}
	\end{equation}
	where $\psi_1=\Gamma(T)\circ(\phi_1)_t$ and $\psi_2=\Gamma(T)\circ(\phi_2)_s$. Here one may notice that $\ol \psi_1=\ol \psi_2$. Since $(L_1)_{t}$ and $ (L_2)_s$ are extended from the base ring $R_{st}$, by Lemma \ref{patchinglemma}, there exists $\omega(T)\in \Aut((R_{st}[T])^m)$ such that (1) $\omega(0)=\text{Id}$ and (2) $\psi_1=\psi_2\circ \omega(T)$. It follows from Quillen's splitting lemma \cite[Theorem 1]{Q} that $\omega(T)$ splits. Therefore, we obtain the following fiber product diagram.
%	$$\begin{tikzcd}
%		0\ar[r] & (L_1)_{t}\ar[rr]\arrow{d} && (A_{st})^m\arrow[d,"\omega(T)"]\arrow[rr,"(\phi_1)_{t}"] && M_{st}\ar[r] \ar[equal]{d} & 0\\
%		0\ar[r] & (L_2)_{s}\ar[rr]&& (A_{st})^m\arrow[rr,"(\phi_2)_s"]  && M_{st}\ar[r] & 0
%	\end{tikzcd}$$
%	Now we consider the following 
	
	$$\begin{tikzcd}
		A^{m} \ar[rrr,dashed] \ar[dd,dashed]\ar[dr,twoheadrightarrow, dashed, "\phi"] &&&  A^{m}_s \ar[dd]\ar[dr, twoheadrightarrow,"\phi_1" ] \\
		& M \ar[rrr,crossing over]  &&& M_s\ar[dd, twoheadrightarrow] \\
		A^{m}_{t} \ar[rr] \ar[dr,twoheadrightarrow,"\phi_2"] & &   A^{m}_{st} \ar[dr,twoheadrightarrow] \ar[r,"\omega(T)"] & A^{m}_{st}\ar[dr,twoheadrightarrow] \\
		& M_{t}\ar[rr] \arrow[uu,crossing over, leftarrow] && M_{st} \ar[r,equal] & M_{st}\ar[from=uu,crossing over]
	\end{tikzcd}$$ 
	Since $\omega(T)$ splits, it follows from the universal property of the fiber product, that there exists a map $\phi:A^{m}\to M$ such that $\phi_s=\phi_1$ and $\phi_{t}=\phi_2$. Since (1) both $\phi_1$ and $\phi_2$ are surjective maps and (2) $\Spec(A)= \Spec(A_{t})\cup \Spec(A_s)$, it follows that $\phi$ is a surjective map. This completes the proof. \qed
		
				Now we are ready to prove the main theorem of this section.
				
	\bt \label{reg}

Let $A$ be a geometric subring of $R[T]$. Let $M$ be an $A$-module. Then $M$ is generated by $e(M)$ many elements.
\et

\proof First we observe some general reductions. Let $\eta$ be the nilradical of $A$ and $A_{\text{red}}=A/\eta$. Since any set of generators of $M/\eta M$ lifts to a set of generators of $M$, without loss of generality we may assume that $A$ is reduced. In particular, this implies that $R$ is reduced. Moreover, using Lemma \ref{l1} one may further assume that $M$ is torsion-free. We give the proof in the following steps.\\
\textbf{Step - 1.} Let $m=e(M)$ and  let ``bar" denote going modulo $\MI_AM$. By \cite[Theorem 0]{P}, we can find $x_1, \ldots, x_m \in \overline{M}$ such that $\overline{M}$ is generated by $x_1, \ldots, x_m$ as a $R$-module. This implies, the map $\phi_1: R^m \twoheadrightarrow \overline{M}$ defined by $\phi_1(e_i) = x_i$ is surjective, where $e_i$'s are the canonical basis elements of $R^m$, $i=1,...,m$. By following the arguments given in [Theorem \ref{ebe}, Step - 1], we can find a non-zero divisor $s \in R$ such that:
\begin{enumerate}[\quad \quad (i)]
	\item $A_s=R_s[T]$ and\label{i5.3}
	\item $M_s$ is a projective module which is extended from $R_s$.\label{ii5.3}
	%	\item $\sigma_1\in \E(A_s\oplus M_s)$ such that $(f,m)\sigma_1=(1,0)$.
\end{enumerate} 
It follows from (\ref{ii5.3}) that any set of generators of $\ol M_s$ can be extended to a set of generators of $M_s$. Therefore, there exists $\omega_i\in M_s$ with $\ol \omega_i=(x_i)_s$ ($i=1,...,m$) such that the map $\phi_1:(A_s)^m\surj M_s$  defined by $\phi_1(e_i)=\omega_i$ is a surjection. Let $L_1=\ker(\phi_1)$. Since $L_1\subseteq A_s^m$, the module $L_1$ is a torsion-free module over $A_s$. Furthermore, repeating the arguments in [Theorem \ref{ebe}, Step - 1] one more time, we may also assume that $(L_1)_s$ is a extended projective $R_s[T]$-module after suitably altering $s$. Note that this alteration of $s$ does not affect any of the favorable properties we obtain in this step.\\ 
\textbf{Step - 2.} 	Let $S=\{1+sr:r\in R\}$. We denote $B=S^{-1}A$ and $L=S^{-1}M$. Let $y_i\in L$ be a lift of $(x_i)_S$, for $i=1,...,m$. Since $s\in R$ is a non-zero divisor, applying Lemma \ref{plgd}, we get a generalized dimension function $d:\Spec(B)\to \mathbb{N}$ such that $d(\p)\le d$ for all $\p\in \Spec(B)$. In particular, this implies $d(\mq)\le \dim(B/\mq) -1$ for all minimal prime ideals $\mq\in \Spec(B)$. Hence we obtain that $m\ge f(M)-1\ge \mu_{\p}(M)+d(\p)$ for all $\p\in {\Spec(B)}$. Therefore, one may apply \cite[Theorem 0]{P} on the triplet  $\{(y_1,...,y_m),S^{-1}(\MI_A M), d\}$, to obtain $z_1,...,z_m\in  L$ such that:
\begin{enumerate}[\quad \quad (a)]
	\item  $z_i=y_i+n_i$, where $n_i\in S^{-1}(\MI_A M)$ and
\item 	$L$ is generated by $z_1,...,z_m$ as a $B$-module.
\end{enumerate}
Then from the construction of $z_i$ it follows that $\ol z_i=x_i$ for all $i=1,...,m$. Let $\phi'_2:B^m\surj L$ be the surjection sending $e_i\to z_i$. Then we observe that $\ol{\phi'_2(e_i)}=x_i$ for all $i=1,...,m$. We define $L'_2:=\ker(\phi'_2)$. Then it follows from $(\ref{ii5.3})$ and from the fact $(L_1)_s$ is a extended projective $R_s[T]$-module, that $(L'_2)_s$ is a stably extended projective module from the base ring $R_{s(1+sR)}$. Let $\mq\in \Spec(B_{s})$ be a minimal prime ideal. Then we note that $A\cap \mq$ is also a minimal prime ideal in $A$ and $(L_{s})_{\mq}=M_{A\cap \mq}$.  Hence we obtain that {\setlength{\abovedisplayskip}{5pt}
 	\setlength{\belowdisplayskip}{5pt}$$\mu_{\mq}({(L'_2)_s})=m-\mu_{\mq}({L_{s}})=m-\mu_{A\cap \mq}(M)\ge \dim(A/A\cap \mq)-1.$$} Now $\dim(B_{s})< \dim(A)$ implies that $\dim(B_{s}/ \mq)\le \dim(A/A\cap \mq)-1$. Therefore, using cancellation result, as stated in \cite[Corollary 1]{P}, one may obtain that $(L'_2)_s$ is an extended projective module from the base ring $R_{s(1+sR)}$. 
 
 There exists an element $t\in S$ and a surjection $\phi_2:(A_t)^m\surj M_t$ such that $\phi_2(e_i)=z_i$ for all $i=1,...,m$ and $(\ker(\phi_2))_s $ is an extended projective module from the base ring $R_{st}$. Let $L_2=\ker(\phi_2)$. Then applying Proposition \ref{pog} one can obtain a surjection $\phi:A^m\surj M $. This concludes the proof. \qed

	\section*{Acknowledgements}
I am grateful to Mrinal Kanti Das for his critical reading of the article and for providing corrections. Above all, I thank him for training me in this subject and always supporting me. Additionally, I extend my gratitude to Md. Ali Zinna for clarifying numerous queries, carefully listening to my lectures regarding this article and offering corrections. This work was financially supported by IISER Kolkata funded Post-doctoral Fellowship.
	
         \bibliographystyle{abbrvurl}
      %\bibliography{Eisenbud-Evans_bbl}

\begin{thebibliography}{10}
			
			\bibitem{BBP}
			S.~Banerjee, C.~Bhaumik, and H.~P. Sarwar.
			\newblock Efficient generation, unimodular element in a geometric subring of a
			polynomial ring. Jan. 2023.
			\newblock \href {https://arxiv.org/abs/2301.11033} {\path{arXiv:2301.11033}}.
			
			\bibitem{SMBAR}
			S.~M. Bhatwadekar and A.~Roy.
			\newblock Stability theorems for overrings of polynomial rings.
			\newblock {\em Inventiones Mathematicae}, 68(1):117--127, 1982.
			\newblock \href {https://doi.org/10.1007/BF01394270}
			{\path{doi:10.1007/BF01394270}}.
			
			\bibitem{SMBB3}
			S.~M. Bhatwadekar and R.~Sridharan.
			\newblock The {E}uler class group of a {N}oetherian ring.
			\newblock {\em Compositio Mathematica}, 122(2):183--222, 2000.
			\newblock \href {https://doi.org/10.1023/a:1001872132498}
			{\path{doi:10.1023/a:1001872132498}}.
			
			\bibitem{MKD}
			M.~K. Das.
			\newblock On a conjecture of {M}urthy.
			\newblock {\em Advances in Mathematics}, 331:326--338, 2018.
			\newblock \href {https://doi.org/10.1016/j.aim.2018.04.012}
			{\path{doi:10.1016/j.aim.2018.04.012}}.
			
			\bibitem{EEC}
			D.~Eisenbud and E.~G. Evans.
			\newblock Three conjectures about modules over polynomial rings.
			\newblock In {\em Lecture Notes in Mathematics}, pages 78--89. Springer Berlin
			Heidelberg, 1973.
			\newblock \href {https://doi.org/10.1007/bfb0068921}
			{\path{doi:10.1007/bfb0068921}}.
			
			\bibitem{EE}
			D.~Eisenbud and E.~G. Evans, Jr.
			\newblock Generating modules efficiently: theorems from algebraic {$K$}-theory.
			\newblock {\em Journal of Algebra}, 27:278--305, 1973.
			\newblock \href {https://doi.org/10.1016/0021-8693(73)90106-3}
			{\path{doi:10.1016/0021-8693(73)90106-3}}.
			
			\bibitem{MPHIL}
			M.~K. Keshari.
			\newblock Euler class group of a {N}oetherian ring.
			\newblock { M-Phil Thesis-2001.} Aug. 2001.
			\newblock \href {https://arxiv.org/abs/1408.2645} {\path{arXiv:1408.2645}}.
			
			\bibitem{NMK}
			N.~M. Kumar.
			\newblock On two conjectures about polynomial rings.
			\newblock {\em Inventiones Mathematicae}, 46(3):225--236, oct 1978.
			\newblock \href {https://doi.org/10.1007/bf01390276}
			{\path{doi:10.1007/bf01390276}}.
			
			\bibitem{m82}
			S.~Mandal.
			\newblock Basic elements and cancellation over {L}aurent polynomial rings.
			\newblock {\em Journal of Algebra}, 79(2):251--257, 1982.
			\newblock \href {https://doi.org/10.1016/0021-8693(82)90301-5}
			{\path{doi:10.1016/0021-8693(82)90301-5}}.
			
			\bibitem{SM1}
			S.~Mandal.
			\newblock On efficient generation of ideals.
			\newblock {\em Inventiones Mathematicae}, 75(1):59--67, 1984.
			\newblock \href {https://doi.org/10.1007/BF01403089}
			{\path{doi:10.1007/BF01403089}}.
			
			\bibitem{MR}
			S.~Mandal and R.~Sridharan.
			\newblock Euler classes and complete intersections.
			\newblock {\em Journal of Mathematics of Kyoto University}, 36(3):453--470,
			1996.
			\newblock \href {https://doi.org/10.1215/kjm/1250518503}
			{\path{doi:10.1215/kjm/1250518503}}.
			
			\bibitem{SM}
			S.~{M}andal (with an appendix~by {M}. {V}.~{N}ori).
			\newblock Homotopy of sections of projective modules.
			\newblock {\em J. Algebraic Geom. 1 (1992)}, (4):639--646, 1992.
			
			\bibitem{E-G}
			J.~McCullough and I.~Peeva.
			\newblock Counterexamples to the {E}isenbud{\textendash}{G}oto regularity
			conjecture.
			\newblock {\em Journal of the American Mathematical Society}, 31(2):473--496,
			nov 2017.
			\newblock \href {https://doi.org/10.1090/jams/891}
			{\path{doi:10.1090/jams/891}}.
			
			\bibitem{MCI}
			M.~P. Murthy.
			\newblock Complete intersections.
			\newblock {\em Conference on Commutative Algebra–1975, pp. 196–211. Queen's
				Papers on Pure and Applied Math., No. 42, Queen's Univ., Kingston, Ont.,
				1975.}
			
			\bibitem{M}
			M.~P. Murthy.
			\newblock Zero cycles and projective modules.
			\newblock {\em Annals of Mathematics. Second Series}, 140(2):405--434, sep 1994.
			\newblock \href {https://doi.org/10.2307/2118605} {\path{doi:10.2307/2118605}}.
			
			\bibitem{P}
			B.~Plumstead.
			\newblock The conjectures of {E}isenbud and {E}vans.
			\newblock {\em American Journal of Mathematics}, 105(6):1417--1433, 1983.
			\newblock \href {https://doi.org/10.2307/2374448} {\path{doi:10.2307/2374448}}.
			
			\bibitem{Q}
			D.~Quillen.
			\newblock Projective modules over polynomial rings.
			\newblock {\em Inventiones Mathematicae}, 36:167--171, 1976.
			\newblock \href {https://doi.org/10.1007/BF01390008}
			{\path{doi:10.1007/BF01390008}}.
			
			\bibitem{RROR}
			R.~Rao.
			\newblock Stability theorems for overrings of polynomial rings. {II}.
			\newblock {\em Journal of Algebra}, 78(2):437--444, 1982.
			\newblock \href {https://doi.org/10.1016/0021-8693(82)90091-6}
			{\path{doi:10.1016/0021-8693(82)90091-6}}.
			
			\bibitem{AMS}
			A.~Sathaye.
			\newblock On the {F}orster-{E}isenbud-{E}vans conjectures.
			\newblock {\em Inventiones Mathematicae}, 46(3):211--224, 1978.
			\newblock \href {https://doi.org/10.1007/BF01390275}
			{\path{doi:10.1007/BF01390275}}.
			
			\bibitem{S}
			J.-P. Serre.
			\newblock Modules projectifs et espaces fibrés à fibre vectorielle.
			\newblock {\em Séminaire Dubreil. Algèbre et théorie des nombres},
			11(2):1--18, 1957-1958.
			\newblock URL: \url{http://eudml.org/doc/111153}.
			
		\end{thebibliography}

		\end{document}